\documentclass[12pt]{article}
\usepackage{amssymb}
\usepackage{amsfonts}
\usepackage{amsmath}
\usepackage{amsbsy}
\newcommand{\be}{\begin{equation}}
\newcommand{\ee}{\end{equation}}
\newcommand{\bea}{\begin{eqnarray}}
\newcommand{\eea}{\end{eqnarray}}
\newcommand{\ba}{\begin{array}}
\newcommand{\ea}{\end{array}}

\newcommand{\bc}{\begin{center}}
\newcommand{\ec}{\end{center}}
\newcommand{\ben}{\begin{enumerate}}
\newcommand{\een}{\end{enumerate}}
\newcommand{\bfi}{\begin{figure}}
\newcommand{\efi}{\end{figure}}

\newcommand{\bq}{\begin{quote}}
\newcommand{\eq}{\end{quote}}
\newcommand{\bqu}{\begin{quotation}}
\newcommand{\equ}{\end{quotation}}
\newenvironment{emphit}{\begin{itemize}}{\end{itemize}}
\newcommand{\bemp}{\begin{emphit}}
\newcommand{\eemp}{\end{emphit}}

\newcommand{\bt}{\begin{tabular}}
\newcommand{\et}{\end{tabular}}

\newtheorem{myth}{Theorem}[section]
\newtheorem{mylem}{Lemma}[section]
\newtheorem{mycor}{Corollary}[section]

\newtheorem{mydef}{Definition}[section]

\newtheorem{myrem}{Remark}[section]

\begin{document}
\date{}
\title{On The Universality Of Central Loops
\footnote{2000 Mathematics Subject Classification. }
}
\author{T\`em\'it\'op\'e Gb\'ol\'ah\`an Ja\'iy\'e\d ol\'a\\Department of
Mathematics,\\
Obafemi Awolowo University, Ile Ife,
Nigeria.\\jaiyeolatemitope@yahoo.com,~tjayeola@oauife.edu.ng}
\maketitle
\begin{abstract}
LC-loops, RC-loops and C-loops are collectively called central
loops. It is shown that an LC(RC)-loop is a left(right) universal
loop. But an LC(RC)-loop is a universal loop if and only if it is a
right(left) universal loop. It is observed that not all RC-loops or
LC-loops or C-loops are universal loops. But if an RC-loop(LC-loop,
C-loop) is universal, then it is a right Bol loop(left Bol loop,
Moufang loop) respectively. If a loop and its right or left isotope
are commutative then the loop is a C-loop if and only if its right
or left isotope is a C-loop. If a C-loop is central square and its
right or left isotope is an alternative central square loop, then
the latter is a C-loop. Necessary and sufficient condition for an
LC-loop(RC-loop) to be a left(right)G-loop is established.
Consequently, necessary and sufficient conditions for an LC-loop,
and an RC-loop to be a G-loop are established. A necessary and
sufficient condition for a C-loop to be a G-loop is established.
\end{abstract}

\section{Introduction}
\paragraph{}
LC-loops, RC-loops and C-loops are loops that respectively satisfy
the identities
\begin{displaymath}
(xx)(yz)=(x(xy))z,\qquad (zy)(xx)=z((yx)x)\qquad \textrm{and}\qquad
x(y(yz))=((xy)y)z.
\end{displaymath}
These three types of loops will be collectively called central
loops. In the theory of loops, central loops are some of the least
studied loops. They have been studied by Phillips and Vojt\v
echovsk\'y \cite{phi, phi1, phi2}, Kinyon et. al. \cite{phi3, phi4,
phi5}, Ramamurthi and Solarin \cite{ram}, Fenyves \cite{fen} and Beg
\cite{phd169, phd170}. The difficulty in studying them is as a
result of the nature of the identities defining them when compared
with other Bol-Moufang identities. It can be noticed that in the
aforementioned LC identity, the two $x$ variables are consecutively
positioned and neither $y$ nor $z$ is between them. A similarly
observation is true in the other two identities(i.e the RC and C
identities). But this observation is not true in the identities
defining Bol loops, Moufang loops and extra loops. Fenyves
\cite{fen} gave three equivalent identities that define LC-loops,
three equivalent identities that define RC-loops and only one
identity that defines C-loops. But recently, Phillips and Vojt\v
echovsk\'y \cite{phi1, phi2} gave four equivalent identities that
define LC-loops and four equivalent identities that define RC-loops.
Three of the four identities given by Phillips and Vojt\v echovsk\'y
are the same as the three already given by Fenyves and their basic
properties are found in \cite{phi, ram, fen, den}. Loops such as Bol
loops, Moufang loops, extra loops are the most popular loops of
Bol-Moufang type whose isotopic invariance have been considered. But
for LC-loops, RC-loops and C-loops, up till this moment, there is no
outstanding result on their isotopic invariance.

The left and right translation maps on a loop $(L,\cdot )$ are the
bijections $L_x:L\to L$ and $R_x:L\to L$, respectively defined as
$yR_x=y\cdot x=yx$ and $yL_x=x\cdot y=xy$ for all $x,y\in L$. The
following subloops of a loop are important for this work.

\paragraph{}
Let $(L, \cdot )$ be a loop. The left nucleus of $L$ is the set
\begin{displaymath}
N_\lambda (L, \cdot )=\{a\in L : ax\cdot y=a\cdot xy~\forall~x, y\in
L\}.
\end{displaymath}
The right nucleus of $L$ is the set
\begin{displaymath}
N_\rho (L, \cdot )=\{a\in L : y\cdot xa=yx\cdot a~\forall~ x, y\in
L\}.
\end{displaymath}
The middle nucleus of $L$ is the set
\begin{displaymath}
N_\mu (L, \cdot )=\{a\in L : ya\cdot x=y\cdot ax~\forall~x, y\in
L\}.
\end{displaymath}
The nucleus of $L$ is the set
\begin{displaymath}
N(L, \cdot )=N_\lambda (L, \cdot )\cap N_\rho (L, \cdot )\cap N_\mu
(L, \cdot ).
\end{displaymath}
The centrum of $L$ is the set
\begin{displaymath}
C(L, \cdot )=\{a\in L : ax=xa~\forall~x\in L\}.
\end{displaymath}
The center of $L$ is the set
\begin{displaymath}
Z(L, \cdot )=N(L, \cdot )\cap C(L, \cdot ).
\end{displaymath}
$L$ is said to be a centrum square loop if~ $x^2\in C(L, \cdot )$
for all $x\in L$. $L$ is said to be a central square loop if~
$x^2\in Z(L, \cdot )$ for all $x\in L$. $L$ is said to be left
alternative if for all $x, y\in L,~ x\cdot xy=x^2y$ and is said to
right alternative if for all $x, y\in L,~ yx\cdot x=yx^2$. Thus, $L$
is said to be alternative if it is both left and right alternative.

The triple $(U, V, W)$ such that $U, V, W\in SYM(L, \cdot )$ is
called an autotopism of $L$ if and only if
\begin{displaymath}
xU\cdot yV=(x\cdot y)W~\forall ~x, y\in L.
\end{displaymath}
$SYM(L, \cdot )$ is called the permutation group of the loop
$(L,\cdot )$. The group of autotopisms of $L$ is denoted by $AUT(L,
\cdot )$. Let $(L, \cdot )$ and $(G, \circ )$ be two distinct loops.
The triple $(U, V, W)~ :~ (L, \cdot )\to (G, \circ )$ such that $U,
V, W : L\to G$ are bijections is called a loop isotopism if and only
if
\begin{displaymath}
xU\circ yV=(x\cdot y)W~\forall ~x, y\in L.
\end{displaymath}
Thus, $L$ and $G$ are called loop isotopes. If the triple $(U, V,
I)~ :~\mathcal{G} =(L, \cdot )\to \mathcal{H}=(L, \circ )$ is an
isotopism, then $\mathcal{H}$ is called a principal isotope of
$\mathcal{G}$. If $U=R_g$ and $V=L_f$, then $\mathcal{H}$ is called
an $f,g$-principal isotope of $\mathcal{G}$. The study of
$f,g$-principal isotopes is important because to every arbitrary
isotope $\mathcal{J}=(G,\ast )$ of a loop $\mathcal{G} =(L, \cdot
)$, there exists an $f,g$-principal isotope $\mathcal{H}=(L, \circ
)$ of $\mathcal{G} =(L, \cdot )$ such that $\mathcal{J}=(G,\ast
)\cong \mathcal{H}=(L, \circ )$. Hence, to verify the isotopic
invariance(or universality) of a loop property, one simply needs to
verify the identity in all $f,g$-principal isotopes. Therefore, a
loop property or identity is isotopic invariant if and only if all
its $f,g$-principal isotopes obey that property. A loop that is
isomorphic to all its loop isotopes is called a G-loop. Thus, a loop
is a G-loop if and only if it is isomorphic to all its
$f,g$-principal isotopes. Some popular G-loops are CC-loops(Goodaire
and Robinson \cite{phd91, phd48}), extra loops, K-loops(Basarab
\cite{phd182, phd137} but not the K-loops of Kiechle \cite{kie}),
VD-loops(Basarab \cite{phd137}) and Buchsteiner loops(Buchsteiner
\cite{phd183}, Basarab \cite{phd46}, Cs\"org\H o et. al.
\cite{phd184}). Some conditions that characterisze a G-loop are
highlighted below.
\paragraph{}
If $H$ is a loop, then the following conditions are equivalent.
\begin{enumerate}
\item $H$ is a G-loop.
\item Every element $x\in H$ is a companion of some right and some
left pseudo-automorphisms of $H$. \qquad Pflugfelder \cite{phd3}
\item There exists a permutation $\theta$ of $H$ such that $(\theta
R_x^{-1},\theta L_y^{-1},\theta )$ is an autotopism of $H$ for all
$x,y$ in $H$. \qquad Chiboka and Solarin \cite{phd45}, Kunen
\cite{phd78}
\end{enumerate}

The study of the $f,g$-principal isotopes of central loops is not
easy because of the nature of the identities defining them as
mentioned earlier. Although its easy to study the $f,g$-principal
isotopes of some other loops of Bol-Moufang type like Bol loops and
Moufang loops. To salvage the situation, the left and right isotopes
of central loops will be studied. According to Nagy and Strambach
\cite{phd88}, the $f,e$-principal isotopes and $e,g$-principal
isotopes of a loop with identity element $e$ are respectively called
its right and left isotopes. The left and right isotopes of a loop
$(\Omega , \cdot)$ shall be represented by $(\Omega ,e,g)$ and
$(\Omega ,f,e)$ respectively.

Most of the expressions that will made in this work are dual in
nature i.e if we write 'statement or word or symbol A(B)[C] implies
a statement or word or symbol
$\mathfrak{A}$($\mathfrak{B}$)[$\mathfrak{C}$]' then we mean that 'A
implies $\mathfrak{A}$', 'B implies $\mathfrak{B}$' and 'C implies
$\mathfrak{C}$'.

In this work, it is shown that an LC(RC)-loop is a left(right)
universal loop. But an LC(RC)-loop is a universal loop if and only
if it is a right(left) universal loop. It is observed that not all
RC-loops or LC-loops or C-loops are universal loops. But if an
RC-loop(LC-loop, C-loop) is universal, then it is a right Bol
loop(left Bol loop, Moufang loop) respectively. If a loop and its
right or left isotope are commutative then the loop is a C-loop if
and only if its right or left isotope is a C-loop. If a C-loop is
central square and its right or left isotope is an alternative
central square loop, then the latter is a C-loop. Necessary and
sufficient condition for an LC-loop(RC-loop) to be a
left(right)G-loop is established. Consequently, necessary and
sufficient conditions for an LC-loop, and an RC-loop to be a G-loop
are established. A necessary and sufficient condition for a C-loop
to be a G-loop is established.

\section{Preliminary}
\begin{mydef}
Let the triple $\alpha =(A,B,C)$ be an isotopism of the groupoid
$(G,\cdot )$ onto a groupoid $(H,\circ )$.
\begin{description}
\item[(a)] If $\alpha =(A,B,B)$, then the triple is
called a left isotopism and the groupoids are called left isotopes.
\item[(b)] If $\alpha =(A,B,A)$, then the triple is
called a right isotopism and the groupoids are called right
isotopes.
\item[(c)] If $\alpha =(A,I,I)$, then the triple is
called a left principal isotopism and the groupoids are called left
principal isotopes.
\item[(d)] If $\alpha =(I,B,I)$, then the triple is
called a right principal isotopism and the groupoids are called
right principal isotopes.
\end{description}
\end{mydef}

\begin{myth}\label{1:1.1}
Let $(G,\cdot )$ and $(H,\circ )$ be two distinct left(right)
isotopic loops with the former having an identity element $e$. For
some $g(f)\in G$, there exists an $e,g$($f,e$)-principal isotope
$(G,\ast )$ of $ (G,\cdot )$ such that $(H,\circ )\cong (G,\ast )$.
\end{myth}

A loop is a left(right) universal "certain" loop if and only if all
its left(right) isotopes are "certain" loops. A loop is a universal
"certain" loop if and only if it is both a left and a right
universal "certain" loop. A loop is called a right
G-loop(G$_\rho$-loop) if and only if it is isomorphic to all its
right loop isotopes. A loop is called a left
G-loop(G$_\lambda$-loop) if and only if it is isomorphic to all its
left loop isotopes. As shown by Wilson \cite{phd90}, a loop is a
G-loop if and only if it is isomorphic to all its right and left
isotopes. Thus, a loop is a G-loop if and only if it is a
G$_\rho$-loop and a G$_\lambda$-loop. Kunen \cite{phd185}
demonstrated the use of G$_\rho$-loops and G$_\lambda$-loops.

\begin{mydef}\label{definition:bijection}(\cite{phd3}, III.3.9~Definition, III.3.10~Definition, III.3.15~Definition)

Let $(L, \cdot )$ be a loop and $U, V, W\in SYM(L, \cdot )$.
\begin{enumerate}
\item If $(U, V, W)\in AUT(L, \cdot )$ for some $V, W$, then $U$ is called autotopic,
\begin{itemize}
\item the set of autotopic bijections in a loop $(L,\cdot )$ is represented by $\Sigma (L,\cdot )$.
\end{itemize}
\item If $(U, V, W)\in AUT(L, \cdot )$ such that $W=U, V=I$, then $U$ is called $\lambda $-regular,
\begin{itemize}
\item the set of all $\lambda $-regular bijections in a loop $(L,\cdot )$ is represented by $\Lambda (L,\cdot )$.
\end{itemize}
\item If $(U, V, W)\in AUT(L, \cdot )$ such that $U=I, W=V$, then $V$ is called $\rho $-regular,
\begin{itemize}
\item the set of all $\rho $-regular bijections in a loop $(L,\cdot )$ is represented by $P(L,\cdot )$.
\end{itemize}
\item If there exists $V\in SYM(L, \cdot )$ such that $xU\cdot y=x\cdot yV$ for all $x, y\in L$,
then $U$ is called $\mu $-regular while $U'=V$ is called its adjoint.
\begin{itemize}
\item The set of all $\mu $-regular bijections in a loop $(L, \cdot )$ is denoted by $\Phi (L, \cdot )$,
while the collection of all adjoints in the loop $(L, \cdot )$ is denoted by $\Phi ^*(L, \cdot )$.
\end{itemize}
\end{enumerate}
\end{mydef}

\begin{myth}\label{lambdarhophi:subgroup}(\cite{phd3}, III.3.11~Theorem, III.3.16~Theorem)

The set $\Lambda (Q,\cdot )\Big(P(Q,\cdot )\Big)\Big[\Phi (Q,\cdot )\Big]$ of all $\lambda $-regular
($\rho $-regular)[$\mu $-regular] bijections of a quasigroup $(Q,\cdot )$ is a subgroup of the group
$\Sigma (Q,\cdot )$ of all autotopic bijections of $(Q,\cdot )$.
\end{myth}

\begin{mycor}\label{lambdarhophi:isomorphism}(\cite{phd3}, III.3.12~Corollary, III.3.16~Theorem)

If two quasigroups $Q$ and $Q'$ are isotopic, then the corresponding groups $\Lambda $
and $\Lambda '$($P$ and $P'$)[$\Phi $ and $\Phi '$]$\{\Phi ^*$ and $\Phi '^*\}$
are isomorphic.
\end{mycor}

Throughout this study, the following notations for translations will
be adopted; $L_x~:~y\mapsto xy$ and $R_x~:~y\mapsto yx$ for a loop
while $L_x'~:~y\mapsto xy$ and $R_x'~:~y\mapsto yx$ for its loop
isotope.

\begin{myth}\label{lc:auto}
A loop $L$ is an LC-loop if and only if $(L_x^2,I,L_x^2)\in AUT(L)$
for all $x\in L$.
\end{myth}
{\bf Proof}\\ Let $L$ be an LC-loop $\Leftrightarrow (x\cdot
xy)z=(xx)(yz)\Leftrightarrow(x\cdot xy)z=x(x\cdot yz)$ by \cite{den}
$\Leftrightarrow(L_x^2,I,L_x^2)\in AUT(L)$ for all $x\in L$.

\begin{myth}\label{rc:auto}
A loop $L$ is an RC-loop if and only if $(I,R_x^2,R_x^2)\in AUT(L)$
forall $x\in L$.
\end{myth}
{\bf Proof}\\ Let $L$ be an RC-loop, then $z(yx\cdot x)=zy\cdot
xx~\Leftrightarrow y(yx\cdot x)=(zy\cdot x)x$ by \cite{den}
$\Leftrightarrow(I,R_x^2,R_x^2)\in AUT(L)~\forall~x\in L$.

\begin{mylem}\label{lcrc:lp}
A loop is an LC(RC)-loop if and only if $L_x^2(R_x^2)$ is
$\lambda(\rho)$-regular i.e $L_x^2(R_x^2)\in\Lambda (L)(P(L))$.
\end{mylem}
{\bf Proof}\\ Using Theorem~\ref{lc:auto}(Theorem~\ref{rc:auto}),
the rest follows from the definition of $\lambda(\rho)$-regular
bijection.

\begin{myth}\label{c:m}
A loop $L$ is a C-loop if and only if $R_x^2$ is $\mu$-regular and
the adjoint of $R_x^2$, denoted by $(R_x^2)^*=L_x^2$ i.e
$R_x^2\in\Phi (L)$ and $L_x^2\in\Phi ^*(L)$.
\end{myth}
{\bf Proof}\\ Let $L$ be a C-loop then $(yx\cdot x)z=y(x\cdot
xz)\Rightarrow yR_x^2\cdot z=y\cdot zL_x^2\Rightarrow R_x^2\in \Phi
(L)$ and $L_x^2\in \Phi ^*(L)$. Conversely: do the reverse of the
above.

\begin{myth}\label{0:12}
Let $G$ be a loop with identity element $e$ and $H$ a quasigroup
such that they are isotopic under the triple $\alpha =(A,B,C)$.
\begin{enumerate}
\item If $C=B$, then $G\overset{A}{\cong}H$ if and only if $eB\in N_\rho(H)$.
\item If $C=A$, then $G\overset{B}{\cong}H$ if and only if $eA\in N_\lambda(H)$.
\end{enumerate}
\end{myth}
{\bf Proof}\\
Here, when $L_x$ and $R_x$ are respectively the left and right
translations of the loop $G$ then the left and right translations of
its quasigroup isotope $H$ are denoted by $L_x'$ and $R_x'$
respectively.

Let $(G,\cdot )$ and $(H,\circ )$ be any two distinct quasigroups.
If $A,B,C : G\rightarrow H$ are permutations, then the following
statements are equivalent :
\begin{itemize}
\item the triple $\alpha=(A,B,C)$ is an isotopism of $G$ upon $H$,
\end{itemize}
\begin{equation}\label{eq:1}
R_{xB}'=A^{-1}R_xC~\forall~x\in G
\end{equation}
\begin{equation}\label{eq:2}
L_{yA}'=B^{-1}L_yC~\forall~y\in G
\end{equation}
\begin{enumerate}
\item When $\alpha=(A,B,B)$, $R_{eB}'=A^{-1}B\Rightarrow
B=AR_{eB}'$. So,
\begin{displaymath}
\alpha=(A,AR_{eB}',AR_{eB}')=(A, A,
A)(I,R_{eB}',R_{eB}'),~\alpha~:~G\rightarrow H.
\end{displaymath}
If $(A, A, A,)~:~G\rightarrow H$ is an isotopism i.e $A$ is an
isomorphism, then $(I,R_{eB}',R_{eB}')~:~H\rightarrow H$ is an
autotopism if and only if $eB\in N_\rho(H)$.
\item When $\alpha=(A,B,A)$, $L_{eA}'=B^{-1}A\Rightarrow
A=BL_{eA}'$. So,
\begin{displaymath}
\alpha=(BL_{eA}',B,BL_{eA}')=(B,B,B)(L_{eA}',I,L_{eA}'),~\alpha~:~G\rightarrow
H.
\end{displaymath}
If $(B, B, B,)~:~G\rightarrow H$ is an isotopism i.e $B$ is an
isomorphism, then $(L_{eA}',I,L_{eA}')~:~H\rightarrow H$ is an
autotopism if and only if $eA\in N_\lambda(H)$.
\end{enumerate}

\section{Main Results}
\begin{myth}\label{iso:lcrc1}
Let $G=(\Omega , \cdot )$ be a loop with identity element $e$ such
that $H=(\Omega , \circ )=(\Omega ,e,g)$ is any of its left
isotopes. $G$ is an LC-loop if and only if $H$ is an LC-loop.
Furthermore, $G$ is a $G_\lambda$-loop if and only if $e\in
N_\rho(H)$.
\end{myth}
{\bf Proof}\\
By Lemma~\ref{lcrc:lp}, $G$ is an LC-loop $\Leftrightarrow L_x^2\in
\Lambda (G)$. From Equation~\ref{eq:2} in the proof of
Theorem~\ref{0:12}, $L_{xg}'=L_x$ for all $g,x\in G$. By
Corollary~\ref{lambdarhophi:isomorphism}, there exists isomorphisms
$\Lambda (G)\to \Lambda (H)$. Thus $L_y'^2\in \Lambda
(H)~\Leftrightarrow~H$ is an LC-loop. The converse is proved in a
similar way.

To prove the last part, we shall use the first part of
Theorem~\ref{0:12}. Thus,  $G$ is a $G_\lambda$-loop if and only if
$e\in N_\rho(H)$.

\begin{mycor}\label{iso:lcrc1.1}
An LC-loop is a left universal loop.
\end{mycor}
{\bf Proof}\\
This follows from Theorem~\ref{iso:lcrc1}.

\begin{mycor}\label{iso:lcrc1.2}
An LC-loop is a
\begin{enumerate}
\item universal loop if and only if it is a right universal loop.
\item G-loop if and only if it is a $G_\rho$-loop and its identity
element is in the right nucleus of all its left isotopes.
\end{enumerate}
 \end{mycor}
{\bf Proof}\\
This follows from Theorem~\ref{iso:lcrc1}.

\begin{myth}\label{iso:lcrc2}
Let $G=(\Omega , \cdot )$ be a loop with identity element $e$ such
that $H=(\Omega , \circ )=(\Omega ,f,e)$ is any of its right
isotopes. $G$ is an RC-loop if and only if $H$ is an RC-loop.
Furthermore, $G$ is a $G_\rho$-loop if and only if $e\in N_\lambda
(H)$.
\end{myth}
{\bf Proof}\\ By Lemma~\ref{lcrc:lp}, $G$ is an RC-loop
$\Leftrightarrow R_x^2\in P(G)$. From Equation~\ref{eq:1} in the
proof of Theorem~\ref{0:12}, $R_{fx}'=R_x$ for all $f,x\in G$. By
Corollary~\ref{lambdarhophi:isomorphism}, there exists isomorphisms
$P(G)\to P(H)$. Thus $R_y'^2\in P(H)~\Leftrightarrow~H$ is an
RC-loop.

To prove the last part, we shall use the second part of
Theorem~\ref{0:12}. Thus, $G$ is a $G_\rho$-loop if and only if
$e\in N_\lambda (H)$.

\begin{mycor}\label{iso:lcrc1.3}
An RC-loop is a left universal loop.
\end{mycor}
{\bf Proof}\\
This follows from Theorem~\ref{iso:lcrc2}.

\begin{mycor}\label{iso:lcrc1.4}
An RC-loop is a
\begin{enumerate}
\item universal loop if and only if it is a left universal loop.
\item G-loop if and only if it is a $G_\lambda$-loop and its identity
element is in the left nucleus of all its right isotopes.
\end{enumerate}
 \end{mycor}
{\bf Proof}\\
This follows from Theorem~\ref{iso:lcrc2}.

\begin{mylem}\label{iso:lcrc1.5}
Not all RC-loops or LC-loops or C-loops are universal loops.
\end{mylem}
{\bf Proof}\\
As shown in Theorem~II.3.8 and Theorem~II.3.9 of \cite{phd3}, a
left(right) inverse property loop is universal if and only if it is
a left(right) Bol loop. Not all RC-loops or LC-loops or C-loops are
right Bol loops or left Bol loops or Moufang loops respectively.
Hence, the proof follows.

\begin{mylem}\label{iso:lcrc1.6}
If an RC-loop(LC-loop, C-loop) is universal, then it is a right Bol
loop(left Bol loop, Moufang loop) respectively.
\end{mylem}
{\bf Proof}\\
This follows from the results in Theorem~II.3.8 and Theorem~II.3.9
of \cite{phd3} that a left(right) inverse property loop is universal
if and only if it is a left(right) Bol loop.

\begin{myrem}
From Lemma~\ref{iso:lcrc1.6}, it can be inferred that an extra loop
is a C-loop and a Moufang loop.
\end{myrem}

\begin{myth}\label{iso:c1}
Let $G=(\Omega , \cdot )$ be a loop with identity element $e$ such
that $H=(\Omega , \circ )=(\Omega ,e,g)$ is any of its left
isotopes. If $G$ is a central square C-loop and $H$ is an
alternative central square loop, then $H$ is a C-loop. Furthermore,
$G$ is a $G_\lambda$-loop if and only if $e\in N_\rho(H)$.
\end{myth}
{\bf Proof}\\ $G$ is a C-loop $\Leftrightarrow R_x^2\in\Phi (G)$ and
$(R_x^2)^*=L_x^2\in\Phi ^*(G)$ for all $x\in G$. From
Equation~\ref{eq:2} in the proof of Theorem~\ref{0:12},
$L_{xg}'=L_x$ for all $g,x\in G$. So using
Corollary~\ref{lambdarhophi:isomorphism}, $R_y'^2 \in\Phi (H)$ and
$(~R_y'^2~)^*=L_y'^2\in\Phi ^*(H) \Leftrightarrow H$ is a C-loop.

The proof of the last part is similar to the way it was proved in
Theorem~\ref{iso:lcrc1}.

\begin{myth}\label{iso:c2}
Let $G=(\Omega , \cdot )$ be a loop with identity element $e$ such
that $H=(\Omega , \circ )=(\Omega ,f,e)$ is any of its right
isotopes. If $G$ is a central square C-loop and $H$ is an
alternative central square loop, then $H$ is a C-loop. Furthermore,
$G$ is a $G_\rho$-loop if and only if $e\in N_\lambda (H)$.
\end{myth}
{\bf Proof}\\ $G$ is a C-loop $\Leftrightarrow R_x^2\in\Phi (G)$ and
$(R_x^2)^*=L_x^2\in\Phi ^*(G)$ for all $x\in G$. From
Equation~\ref{eq:1} in the proof of Theorem~\ref{0:12},
$R_{fx}'=R_x$ for all $f,x\in G$. So using
Corollary~\ref{lambdarhophi:isomorphism}, $R_y'^2 \in\Phi (H)$ and
$(~R_y'^2~)^*=L_y'^2\in\Phi ^*(H) \Leftrightarrow H$ is a C-loop.

The proof of the last part is similar to the way it was proved in
Theorem~\ref{iso:lcrc2}.

\begin{myth}\label{iso:cc1}
Let $G=(\Omega , \cdot )$ be a commutative loop with identity
element $e$ such that $H=(\Omega , \circ )=(\Omega ,e,g)$ is any of
its commutative left isotopes. $G$ is a C-loop if and only if $H$ is
a C-loop. Furthermore, $G$ is a $G_\lambda$-loop if and only if
$e\in N_\rho(H)$.
\end{myth}
{\bf Proof}\\ The proof is similar to that of Theorem~\ref{iso:c1}.

\begin{myth}\label{iso:cc2}
Let $G=(\Omega , \cdot )$ be a commutative loop with identity
element $e$ such that $H=(\Omega , \circ )=(\Omega ,f,e)$ is any of
its commutative right isotopes. $G$ is a C-loop if and only if $H$
is a C-loop. Furthermore, $G$ is a $G_\rho$-loop if and only if
$e\in N_\lambda (H)$.

\end{myth}
{\bf Proof}\\ The proof is similar to that of Theorem~\ref{iso:c1}.

\begin{mylem}\label{}
A central loop is a G-loop if and only if its identity element
belongs to the intersection of the left cosets formed by its
nucleus.
\end{mylem}
{\bf Proof}\\
Let $G$ be a central loop. Let its left isotopes be denoted by
$H_i,~i\in \Pi$ and let its right isotopes be denoted by $H_j,~j\in
\Delta$. By Theorem~\ref{0:12}, $G$ is a $G_\lambda$-loop if and
only if $e\in N_\rho(H_i)$ for all $i\in \Pi$. Also, $G$ is a
$G_\rho$-loop if and only if $e\in N_\lambda (H_j)$ for all $j\in
\Delta$. A loop is a G-loop if and only if it is a G$_\rho$-loop and
a G$_\lambda$-loop. So $G$ is a G-loop if and only if $e\in
\bigcap_{{}_{i\in \Pi ,j\in \Delta}}\Big(N_\rho(H_i)\cap N_\lambda
(H_j)\Big)$.

$N_\rho(H_i)$ and $N_\lambda (H_j)$ are subgroups for all $i\in \Pi$
and $j\in \Delta$. Let $G=(\Omega , \cdot )$ so that $H_i=(\Omega ,
\circ_i )=(\Omega ,e,g_i),~i\in \Pi$ and $H_j=(\Omega , \circ_j
)=(\Omega ,f_j,e),~j\in \Delta$. Recall that
$N_\rho(G)\overset{L_{g_i}}{\cong}N_\rho(H_i)$ for all $i\in \Pi$
i.e $N_\rho(G)\cong N_\rho(H_i)$(III.2.6~Theorem \cite{phd3}) under
the mapping $L_{g_i}$. Similarly, $N_\lambda
(G)\overset{R_{f_j}}{\cong}N_\lambda (H_j)$ for all $j\in \Delta$
i.e $N_\lambda (G)\cong N_\lambda (H_j)$(III.2.6~Theorem
\cite{phd3}) under the mapping $R_{f_j}$. Therefore, $G$ is a G-loop
if and only if $e\in \bigcap_{{}_{i\in \Pi ,j\in
\Delta}}\Big(g_iN_\rho(G)\cap N_\lambda (G)f_i\Big)$.

In the case of $G$ been a C-loop, $N(G)=N_\rho(G)=N_\lambda (G)$,
$N(G)\vartriangleleft G$ and $[G~:~N(G)]\ne 2$. So, $G$ is a G-loop
if and only if $e\in \bigcap_{{}_{g\in \{g_i,f_j\}}}^{{}^{\{i\in \Pi
,j\in \Delta\}}}gN$.

\begin{mycor}
Every alternative central square left(right) isotope $G$ of a Cayley
loop or RA-loop or $M_{16}(Q_8)\times E\times A$ where $E$ is an
elementary abelian~2-group, $A$ is an abelian group(all of whose
elements have finite odd order) and $M_{16}(Q_8)$ is a Cayley loop,
is a C-loop.
\end{mycor}
{\bf Proof}\\
From \cite{goo1}, the Cayley loop and ring alternative
loops(RA-loops) are all central square. Hence, by
Theorem~\ref{iso:c1} and Theorem~\ref{iso:c2}, the claim that $G$ is
a C-loop follows.

\end{document}